\documentclass[12pt]{amsart}
\usepackage[dvips]{graphicx,color}   
\usepackage{amsmath}
\usepackage{amscd}
\usepackage{amsfonts,latexsym,amssymb}

\DeclareFontFamily{OT1}{rsfs}{}
\DeclareFontShape{OT1}{rsfs}{n}{it}{<->rsfs10}{}
\DeclareMathAlphabet{\curly}{OT1}{rsfs}{n}{it}

\newcommand{\CrC}{{\curly C}}

\newcommand{\CrN}{{\curly N}}

\newcommand{\CrP}{{\curly P}}

\newcommand{\Z}{\mathbb{Z}}
\newcommand{\C}{\mathbb{C}}

\newcommand{\N}{\mathbb{N}}
\newcommand{\R}{\mathbb{R}}

\newcommand{\CP}{\mathbb{CP}}

\newtheorem{proposition}{Proposition}[section]
\newtheorem{theorem}[proposition]{Theorem}
\newtheorem{definition}[proposition]{Definition}
\newtheorem{lemma}[proposition]{Lemma}

\newtheorem{corollary}[proposition]{Corollary}

\def\>#1{{\bf #1}}             

\def\dim{\hbox{{\rm dim}}}

\def\Gr{\hbox{{\rm Gr}}}
\def\Ker{\hbox{{\rm Ker}}}
\def\length{\hbox{{\rm length}}}
\begin{document}

{\mbox{   }}

\vskip 2cm

\centerline{\Large\bf LEFSCHETZ TYPE PENCILS}

\bigskip

\centerline{\Large\bf ON CONTACT MANIFOLDS}

\medskip

\vskip 2cm

\centerline{F. Presas}

\vskip 1.5cm

\medskip
\centerline{ {\it Departamento de \'Algebra, Universidad 
Complutense de Madrid, 28040 Madrid, Spain.}  }

\vskip 2cm

\centerline{\sc Abstract}
{\small We define the concept of Lefschetz contact pencil and we show the existence
of such structures on any contact manifold. The main idea of the proof is a 
generalization of the Donaldson arguments used in the symplectic case. 
We will analyze some of the applications of such existence theorem for the topology of
approximately holomorphic contact submanifolds. In particular, we study the topological
relationship between the contact submanifods constructed in \cite{IMP99}.}


\vskip 1.5cm

\newpage

\section{Introduction.} \label{introduction}

S. Donaldson in \cite{Do96} has adapted the concept of very ample bundle to the
symplectic setting. Following this idea some results of complex projective geometry
has been generalized to symplectic geometry, as Bertini's theorem and Lefschetz's
hyperplane theorem in \cite{Do96}, conectedness of the space of ``good'' sections of a very
ample bundle in \cite{Au97}, divisors on projective fibrations in \cite{Pa98},
special position theorems in \cite{Pa99}, existence of Lefschetz pencils in 
\cite{Do99}, existence of branched coverings of symplectic 4-manifolds over $\CP^2$
and associated invariants \cite{Au99, Au99b}, and Kodaira's embeddings theorems and 
symplectic determinantal submanifolds in \cite{MPS99}. These ideas have opened a new
insight in symplectic geometry allowing to understand symplectic manifolds through 
the study of the linear systems associated to a ``very ample'' vector bundle.

In \cite{IMP99} the idea of \cite{Do96} of working with approximately $J$-holomorphic sections
was partially traslated to the contact case.
With these sections and a generalization of the local estimated transversality
result of \cite{Do96}, which is the key of the symplectic approach, the ideas of 
\cite{Do96,Au97} were translated word by word to the contact setting. The only important
loss was the isotopy results which have been developed in the symplectic theory starting with
the ideas of D. Auroux in \cite{Au97}. In this paper we show how to develop a contact 
geometry of linear systems analogous to the symplectic case. We also show how to recover 
partially the isotopy in the theory.

We will prove a theorem analogue to that of \cite{Do99}. In fact, we will show the 
existence of a certain class of pencils on a contact manifold. The main tool in the proof 
will be a generalization of the local transversality theorem proved in \cite{Do99}. As in
the symplectic case this result is not strictly necessary for the proof, but it has interest
in its own. It could be used to simplify the constructions of other linear systems,
as for instance, the one in \cite{MPS99} in the symplectic case.

A compatible chart in a contact manifold $(C,D)$ at a point $x$ will be a chart 
$\phi: U_x \subset C \to \C^n\times
\R$, where $U_x$ is neighborhood of $x$, verifying $\phi(x)=(0,0)$, $(\phi_*)(D(x))=\C^n\times \{ 0 \}$ and moreover verifying that the 
presymplectic form $(\phi_*)d \theta(x)$, when restricted to $\C^n\times\{ 0 \}$, is a 
positive form of type $(1,1)$ at the origin of coordinates. If the contact manifold is exact 
we will impose 
also that $\langle \phi_*(R)(x)), \frac{\partial}{\partial s} \rangle>0$ and say that 
the chart is 
oriented compatible, where $R$ is the Reeb vector field and $s$ is the real coordinate.
\begin{definition} \label{pencil_def}
A (oriented) contact pencil on a closed (exact) contact manifold $C$ consists of the 
following data:
\begin{enumerate}
\item a codimension $4$ contact submanifold $A\subset C$,
\item a finite set of smooth contact curves $\Delta= \bigcup_{i\in I} \gamma_i \subset C-A$,
\item  a smooth map $f: V-A \to \CP^1$, whose restriction to the complementary
of the set $\Delta$ is a submersion, satisfying also that $P_\Delta=f(\Delta)$ is a set
of locally smooth curves with transversal selfintersections.
\end{enumerate}
Also the data have to admit the following standard local models:
\begin{itemize}
\item At any point $a\in A$, there are (oriented) compatible coordinates 
$(z_1, \ldots,z_n, s)\in \C^n \times \R$ such that $A$ is locally given by $\{z_1,z_2=0\}$.
And the function $f$ has the expression $f(z_1, \ldots, z_n, s)=
\frac{z_1}{z_2}\in \CP^1$ near $a$.
\item At a point $b_i \in \gamma_i$ there are (oriented) compatible
coordinates in which $f$ is written as $f(b_i)+\varphi(s)+z_1^2+\ldots+z_n^2$,
where $\varphi:\R \to \C$ verifies $\varphi(0)=0$ and $\varphi'(0)\neq 0$.
\end{itemize}
\end{definition}
It is clear from the local model that the counterimage $f^{-1}(p)$ is a subset of $C-A$, 
whose closure in $C$ is smooth at $f^{-1}(p)\bigcap A$. Abusing language, we will
call fiber over
$b$ to the closure of the counterimage. It is a smooth
submanifold if $b$ is a regular value of $f$. In other case we will have one or two 
singularities locally modelled by:
\begin{eqnarray} 
\varphi(s)+z_1^2+\cdots +z_n^2. \label{model_sin}
\end{eqnarray}
In case $\dim ~ C=3$, the smooth fibers will be (oriented) links on $C$. The link 
operation that is performed when the image crosses a circle $\gamma_{i}$ of the sphere
$\CP^1$ looks, after general projection to a plane, as

\begin{figure}[htbp]
\centerline{
\includegraphics[
width=12.5cm,       
]
{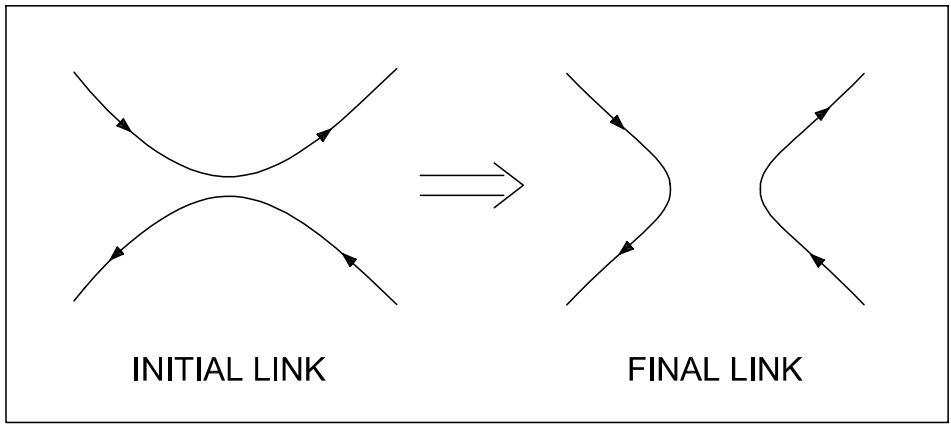}            
}
\caption{\em Link operation}

\end{figure}

The main result of this article is
\begin{theorem} \label{main_th}
Given a contact closed manifold $(C,D)$ (resp. exact) and $\alpha\in H_{2n-1}(C,\R)$
which is reduction of an integer class,
there exists a  contact pencil on $C$ (resp. oriented) whose fibers are contact submanifolds,
homologous to $\alpha$.
\end{theorem}

With this result at hand, it is easy to understand the possibility
of finding general isotopic constructions for the contact submanifolds constructed 
in \cite{IMP99}. Once fixed a compatible complex structure in the distribution, in Section 
\ref{topology} we will show how to define a sequence of contact
fibrations $f_k$, satisfying  $f_k(0)\simeq N_k$ and $f_k(\infty)\simeq N_k'$,
where $N_k$ and $N_k'$ are sequences of codimension $2$ contact submanifolds constructed 
with the method developed in \cite{IMP99}. In order to assure the isotopy between $N_k$ and 
$N_k'$ we have only
to construct a path between $0$ and $\infty$ in $\CP^1$ which does not intersect
$\Delta$, because in this case $f$ restricted to the path is surjective. But this is only 
possible if $0$ and 
$\infty$ are in the same connected component
of $\CP^1-\Delta$, which is not true in general (contrarely to the symplectic
case where $\Delta$ consists of isolated points). We will study in Subsection 
\ref{previous} the topological relationship between the counterimages of the points of a 
path crossing $\Delta$. This will prove that
$$ H_i(N_k)=H_i(N_k'),$$
$$ \pi_i(N_k)=\pi_i(N_k'),$$
for $k=0, \ldots, n-2$. Moreover we will show how to make the result independent of the chosen
complex structure and of the contact form, thus providing a comprobation of the
``contact Lefschetz hyperplane theorem''
proved in \cite{IMP99}. The result brings us back to the Lefschetz original ideas to prove
the hyperplane theorem (see \cite{Le24}). It would be interesting to complete the proof in
Lefschetz's way, since this would be a more geometrical proof than the one contained in 
\cite{IMP99}, the latter
being an adaptation of the modern Morse-theoretic argument due to Bott, Andreotti and Frankel
(and in the symplectic case to Donaldson).

In Section \ref{prelim} we will give the basic results in contact geometry needed to 
develop the proof.
In Section \ref{proof} we will give the proof of the exact case assuming
transversality results of approximately holomorphic contact geometry. 
Afterwards, in Sections \ref{asymp} and \ref{local_res}, we
develop the local approximately holomorphic techniques needed to achieve this transversality.
In Section \ref{topology} we will study the relationship between approximately
holomorphic sections transverse to {\bf 0}.
Finally, we will adapt all the study to the non-exact case in Section \ref{double}.

{\bf Acknowledments:} \newline
I am very grateful to S. Donaldson by his kindness passing me the preprint \cite{Do99}
which is absolutely fundamental to develop this work and also by his useful observations
about some of the parts of this paper. Also I want to mention 
D. Auroux and R. Paoletti for passing me their useful preprints. Finally
we want to thank to the members of the GESTA \footnote{Geometr\'{\i}a simpl\'ectica con
t\'ecnicas algebraicas.} seminar in Madrid their support and interest through the
elaboration of this work. I want to mention especially Vicente Mu\~noz by his
useful commentaries to the previous versions of this paper.

\section{Definitions and results.} \label{prelim}
We will assume along the proofs that $(C,D)$ is an exact contact manifold, where we have
fixed a global contact form $\theta$. In Section \ref{double} we will precise the changes 
needed to extend the results to the non-exact setting. We start recalling with
the basic definitions and results of \cite{IMP99}.

\subsection{Basic concepts}
A complex structure on $C$ is a complex structure $J$ defined on $D$, interpreted as a
symplectic vector bundle, compatible with the symplectic form $d\theta$. Recall that 
the distribution $D$ is contact, if and only if the restriction of $d\theta$ to $D$ is 
symplectic. The contact $1$-form defines a vector field $R$ by the condition
$$ i_R\theta=1, ~~~ i_R d\theta=0, $$
It is called the Reeb vector field. Given a complex structure we obtain a metric for 
the manifold
$$ g_J(u,v)= d\theta(u,Jv)+ \theta\otimes \theta. $$
This metric is called a contact metric. It depends on the fixed contact form.
The $k$-rescaled contact metric will be defined as $g_{J,k}= kg_J$. If there is not risk
of confusion we will denote simply $g_k$, supposing $J$ fixed along the proofs.
(Observe that the $k$-rescaled contact metric $g_k$ is not the contact metric associated to
$k^{\alpha}\theta$, for $\alpha\in \R$) 

Now we give some definitions to control contact structures in $\R^{2n+1}$.
\begin{definition}
The maximum angle between two subspaces $U,V\in \Gr_{\R}(r,n)$ is defined as:
$$ \angle_M(U,V)= \max_{u\in U} \angle(u,V). $$
\end{definition}
This angle defines a distance in the topological space $\Gr_{\R}(r,n)$ (for details
see \cite{MPS99}).
\begin{definition}
Let $\theta_k$ be a sequence of contact forms in $\R^{2n+1}$, with associated distributions
$D_k=\Ker ~ \theta_k$. The sequence is called $c$-asymptotically flat in the set 
$U\subset \R^{2n+1} $ if
$$ \angle_M(D_k(0), D_k(x))\leq ck^{-1/2}, ~~ \hbox{{\rm for all }} x\in U. $$
The sequence is called asymptotically flat if there exists some $c>0$ for which it is 
$c$-asymptotically flat.
\end{definition}
The standard contact structure in $\R^{2n+1}$ is defined as $\theta_0=ds+\sum_{j=1}^n x_jdy_j$, 
where $(x_j,y_j,s)\in \R^{2n+1}$. The $k$-rescaled contact metric is the contact metric
associated to $\theta_{k^{1/2}}=k^{1/2}ds+\sum x_jdy_j$, which is obtained from $\theta_0$
scaling the coordinates by a factor $k^{1/2}$. It is obvious from the definition that
$\theta_{k^{1/2}}$ is a sequence of asymptotically flat contact forms on any bounded set
of $\R^{2n+1}$.

If a contact distribution $D$ is defined at the origin by the horizontal 
subspace $\R^{2n}\times \{ 0 \}\subset \R^{2n+1}$, we can define at a neighborhood  of 
this point a canonical complex structure in $D$ by means of vertical projection of the canonical
one defined on $\R^{2n}$. This complex structure will be denoted $J_0$. In the case
of the standard contact structure we can extend $J_0$ over all $\R^{2n+1}$.

Following \cite{IMP99} we can define $\partial$ and $\bar{\partial}$ operators in
any kind of function, morphism or section defined on a contact manifold, restricting 
ourselves to the contact distribution, or equivalently projecting along the Reeb direction. 

\subsection{Approximately holomorphic geometry}
We will talk about uniform constant, polynomial, when this constant, polynomial, etc.
does not depend on the chosen point $x\in C$, nor in the integer $k$ appearing in the context.
(However these constants can depend on the modulus of the given sections, on the nature of its
derivatives, on the step of a recurrent reasoning, etc. but always independently of $k$).

We will use from \cite{IMP99} the following
\begin{lemma} \label{quasi}
Given $(C,\theta)$ a closed contact manifold and $J$ a compatible complex structure on $\Ker ~
\theta$, there exists a uniform constant $c>0$ and a contact
Darboux chart $\psi: (C,\theta) \to (\R^{2n+1}, \theta_0)$ satisfying
$$ \frac{1}{2}g(v,w) \leq \langle (\psi_*)_x v,(\psi_*)_xw \rangle \leq 2g(v,w),
\forall x\in B_g(x,c), v,w\in T_xC. $$
This implies that $|\nabla^r \psi|=O(1)$ and $|\nabla^r \psi^{-1}|=O(1)$, for $r=1,2,3$.
Moreover $|\bar{\partial} \psi(y)|\leq c'd(x,y)$, for a uniform constant $c'$.
\end{lemma}
A $k^{1/2}$-Darboux chart is a chart $\phi_k:B(x, \epsilon) \to \R^{2n+1}$ such
that $(\phi_k)_* \theta =\theta_{k^{1/2}}$. The following result is a direct corollary of 
Lemma \ref{quasi}:
\begin{corollary} \label{approx}
Given $(C,\theta)$ a closed contact manifold and $J$ a compatible complex structure on $\Ker
~ \theta$, there exists a uniform constant $c>0$ such that there exists a contact
chart $\psi: (C,\theta) \to (\R^{2n+1}, \theta_{k^{1/2}})$ satisfying
$$ \frac{1}{2}d_{g_k}(v,w) \leq \langle (\psi_*)_x v,(\psi_*)_xw \rangle \leq 2d_{g_k}(v,w),
\forall x\in B_{g_k}(x,c), v,w\in T_xC. $$
This implies that $|\nabla^r \psi|=O(1)$ and $|\nabla^r \psi^{-1}|=O(1)$, for $r=1,2,3$.
Moreover $|\nabla^r \bar{\partial} \psi(y)|\leq c'k^{-1/2}$, for $r=0,1,2$, for a uniform 
constant $c'$.
\end{corollary}
These results allows to trivialize locally contact manifolds in an approximately holomorphic
way. The analogous notion of the $c$-bounds of \cite{Do99} in the contact case is the 
following
\begin{definition}
A sequence of sections $s_k$ of hermitian bundles $E_k$ over the contact manifold 
$(C,\theta)$ has mixed $C^r$-bounds $(c_D, c_R)$ at the point $x\in C$ if it satisfies
\begin{eqnarray*}
& & |s_k(x)|<c_D, \\
& & |\nabla_D^j s_k(x)|<c_D, ~~ \forall j=1, \ldots, r. \\
& & |\nabla^j s_k(x)|< c_R, ~~ \forall j=1, \ldots, r. \\
& & |\nabla^j \bar{\partial} s_k(x)| <c_R k^{-1/2}, ~~ \forall j=0, \ldots, r-1.
\end{eqnarray*}
The sequence has uniform mixed $C^r$-bounds $(c_D, c_R)$ if it satisfies these bounds at every
point.
\end{definition}
There we denote by $\nabla_D$ the restriction of the operator $\nabla$ to the subspace $D$. 
The metric used in the manifold $C$ to measure the norms in the precedent definition is
the rescaled contact metric $g_k$.
Recall that if $s_k^j$ has $(c_D^j, c_R^j)$ mixed $C^r$-bounds ($j=1,2$), then 
$s_k^1+s_k^2$ has $(c_D^1+c_D^2, c_R^1+c_R^2)$ mixed $C^r$-bounds.

The following results are used to trivialize bundles over contact manifolds.
\begin{definition} 
A sequence of sections $s_k$ of bundles $E_k$ has mixed Gaussian decay in
$C^r$-norm away from a point $x\in C$ if there exist a uniform polynomial $P$ and
uniform constants $\lambda > 0$, $c_a$ such that 
for all $y\in X$, the sequence $s_k(y)$ has mixed $C^r$-bounds 
$$(P(d_k(x,y)) \exp (-\lambda d_k (x,y)^2), c_aP(d_k(x,y)) \exp (-\lambda d_k (x,y)^2)).$$
\end{definition}
Recall that the uniform constant $c_a$ is necessary because we will check afterwards that 
for different sequences of sections this constant cannot be fixed.

The prequantizable line bundle $L$ over an exact contact manifold $C$ is defined as
the complex line bundle with connection such that $curv(L)=id\theta$.
\begin{lemma}[Lemma 5 from \cite{IMP99}] 
Let $(C, \theta )$ be a closed contact
manifold.  There exists a uniform  constant
$c_s>0$, such that given any point $x\in C$, there exists a sequence
of sections $\sigma_{k,x}$ of $L^{\otimes 
k}$ satisfying $|\sigma_{k,x}| \geq c_s$ at every $y$ in a ball of
$g_k$-radius $10$ centered at $x$ and the
sections $\sigma_{k,x}$ have uniform mixed Gaussian decay away from $x$ in $C^3$-norm
(in this case $c_a=1$).
\end{lemma}

\subsection{Transversality results.}
Following Donaldson \cite{Do99} and Auroux \cite{Au99} we set up the following definitions.
A linear application $f: \R^n \to \R^r$ is $\eta$-transverse if it has a right inverse
$\nu:\R^r \to \R^n$ such that $|\nu|\leq \eta^{-1}$. In the non-linear case we
will say that $f:U\subset \R^n \to \R^r$ is $\eta$-tranverse to $y\in \R^r$ over $U$ if
$\forall x \in U$, such that $|f(x)-y|<\eta$, then $df$ is $\eta$-transverse. 
Recall that this is an open condition. In fact, if $|f-g|_{C^1, U}<\epsilon/10$ and $f$
is $\epsilon$-transverse to $y$ over $U$ then $g$ is, say, $\epsilon/2$-transverse to
$y$ in $U$.

The definition of transversality to {\bf 0} for sections of hermitian bundles over 
riemannian manifolds is totally analogous. In the case of contact manifolds we
have to strengthen the conditions.
\begin{definition}
A section $s_k$ of the hermitian vector bundle $E_k$ over the contact manifold $(C, \theta)$
is $\eta$-transverse to {\bf 0} on $U\subset C$ if for all $x\in U$ such that $|s_k(x)|
<\eta$ then $\nabla_D s_k(x)$ is $\eta$-transverse to {\bf 0} 
(with respect to the $g_k$-metric in $C$).
\end{definition}
From the discussion of \cite{IMP99} it follows that a sequence of sections $s_k$ of the
bundles $E_k$ over the contact manifold $(C,\theta)$, which has uniform $(c_D, c_R)$ bounds 
and which is $\eta$-transverse to {\bf 0}, has as zero set a contact submanifold, for $k$
large enough. The precise result of \cite{IMP99} is
\begin{theorem} \label{trans_1}
Given a closed exact contact manifold $(C, \theta)$. Let $\epsilon>0$ and let $s_k$ be a 
sequence of sections of the bundles $E\otimes L^{\otimes 
k}$, for a fixed hermitian bundle $E$, with 
uniform mixed $C^r$-bounds $(c_D, c_R)$. Then there exists a real number $\eta>0$ (depending on
$\epsilon$, $c_D$ and $c_R$), and a sequence $\sigma_k$ such that:
\begin{enumerate}
\item $\sigma_k-s_k$ has mixed $C^r$-bounds $(\epsilon,c_R')$.
\item $\sigma_k$ is $\eta$-transverse to {\bf 0}.
\end{enumerate}
\end{theorem}
In \cite{IMP99} the proof is developed for sequences with mixed $C^2$-bounds, but there is
not any problem in generalizing it to the mixed $C^r$-bounds case. We need in this
article mixed $C^3$ bounds.
Given a section $s_k=(s_k^0,s_k^1)$ of the bundle $\C^2\otimes S\otimes L^{\otimes k}$, for 
a fixed hermitian line bundle $S$, whose zero set is $Z(s_k)$, we denote
$F^{s_k}=\frac{s_k^1}{s_k^0}: C- Z(s_k) \to \CP^1$
the projectivization of the section. The holomorphic part of the differential of this
application will be denoted by $\partial F^{s_k}$.
Now we state a generalization of Theorem \ref{trans_1} which will be proved in Section 
\ref{asymp}:
\begin{theorem} \label{trans_2}
Given a closed exact contact manifold $(C, \theta)$.
Let $\epsilon>0$ and let $s_k=(s_k^0,s_k^1)$ be a sequence of sections of the bundles 
$\C^2\otimes S \otimes L^{\otimes k}$ with
mixed $C^3$-bounds $(c_D, c_R)$. Suppose that $s_k^0$ and $s_k$ are both transverse to {\bf 0}.
Then there exists $\eta>0$ (depending on
$\epsilon$, $c_D$ and $c_R$), and a sequence $\sigma_k$ satisfying
\begin{enumerate}
\item $\sigma_k-s_k$ has mixed $C^3$-bounds $(\epsilon,c_R')$,
\item $\sigma_k^0$ and $\sigma_k^0\oplus \sigma_k^1$ are $\eta$-transverse to {\bf 0},
\item $\partial F_k^{\sigma}$ is $\eta$-transverse to {\bf 0} away from $Z(\sigma_k^0)$.
\end{enumerate}
\end{theorem}
The techniques used in Section \ref{asymp} improve slightly the ones of \cite{IMP99}
and thus could allow a simpler proof of Theorem \ref{trans_1} avoiding some of the 
complications of the globalization process in that article.

\section{Proof of the main result.} \label{proof}
Take a complex line bundle $S$ with connection $\nabla$ satisfying that 
$curv(\nabla)=PD(\alpha)$.
This is possible since $\alpha$ is an integer class.
Starting with any sequence of sections $s_k'$ of $\C^2\otimes S\otimes 
L^{\otimes k}$ with mixed $C^3$-bounds, we can perturb it using Theorems
\ref{trans_1} and \ref{trans_2} to achieve a sequence of sections $s_k$ verifying
properties 2 and 3 of Theorem \ref{trans_2} and with mixed $C^3$-bounds $(c_D, c_R)$. 
We consider this sequence as starting datum and
will use it to construct the oriented contact pencil.

From the $\eta$-transversality of $s_k$ the zero set $A=Z(s_k)$ is a codimension $4$ 
contact 
manifold where $F^{s_k}$ is not well defined. We will write $F$ instead of $F^{s_k}$ 
whenever it causes no confusion. Now we will study, as in \cite{Do99}, the shape of the ``bad
set'' $\Gamma=\{ x\in C: |\partial F|\leq |\bar{\partial} F| \}$. From \cite{IMP99}
we know that if we prove that $\Gamma= \{x\in C: d_DF=0 \}$ we will have obtained
that the fibres of $F$ are contact at all smooth points. This will be proved in several
steps. 
\begin{lemma} \label{far}
There is a constant $\xi>0$, depending only on $c_D$, $c_R$ and $\epsilon$, such that if $k$ is
large enough then $|s_0|\geq \xi$ on $\Gamma$.
\end{lemma}

{\bf Proof:}
The proof is analogous to that of Lemma 7 in \cite{Do99}. \hfill $\Box$

Define $\Delta$ as the set of critical points of $F$.
The connected components of $\Delta$ form a discrete set of smooth curves by the 
transversality condition imposed to $\partial F$. Also we can assure that this set of 
components is finite
because it is contained in $\Gamma$, which by Lemma \ref{far} is contained in the 
complementary of a $\gamma$-neighborhood of $W_{\infty}=Z(s_k^0)$, for $\gamma>0$
a uniform constant small enough. We define:
$$ \Omega_{\xi}= \{ p\in C, ~ |s_0(p)|> \xi/2 \}. $$

The following step, adapting again Donaldson's argument, is to estimate the shape of
the set $\Delta$. This is the content of the following:
\begin{proposition} \label{disjoint}
There is a uniform constant $\rho_0>0$, such that the $\rho_0$-neighborhoods
of each connected component $\gamma_i$ of $\Delta$ are disjoint and are contained
in $\Omega_{\xi}$. Moreover for any $\rho<\rho_0$, for $k$ large enough (only depending
on $\rho$), the set $\Gamma$ is contained in a $\rho$-neighborhood of the set $\Delta$.
\end{proposition}
The proof of this result is absolutely analogous to the proof of Proposition 9 in 
\cite{Do99} and depends strongly on Lemma 8 in that paper. We refer the reader to
\cite{Do99} for the argument.

Finally, we need to perturb the sequence $\partial F^{s_k}$ in arbitrarily small neighborhoods
of $A$ and $\Delta$ to achieve the local models required in Definition \ref{pencil_def}.
The perturbation required in $\Delta$ needs a careful analysis, but again the situation
in $A$ is a straightforward generalization of \cite{Do99}. We need only to define
$$ LD_x= \nabla s_0^k \oplus \nabla s_1^k: TC_x \to L_x^{\otimes k} \oplus
L_x^{\otimes k}. $$
With this notation the result we need in our case is
\begin{lemma} \label{positi}
For a point $x\in A$, $F$ can be represented in the standard model of Definition 
\ref{pencil_def} at $x$ if and only if $TA_x\bigcap D$ is a symplectic subspace and the 
restriction of $d\theta$ to the symplectic orthogonal $CA_x$ (in $D$) of $TA_x\bigcap D$ is a
positive form of type $(1,1)$ with respect to the complex structure on $CA_x$ induced
by $LD_x$.
\end{lemma}
We do not provide a proof of this Lemma, which follows word by word the proof of Lemma 11 
in \cite{Do99}.
Thanks to the $(c_D',c_R')$ mixed $C^3$-bounds and the transversality of the sequence
$s_k$, it is easy to check that a small $C^3$-perturbation of the sequence satisfies
the hypothesis of Lemma \ref{positi}, thus completing the study in the neighborhood of $A$.

Now we study the map $F$ near $\Delta$. Again Donaldson's ideas work in this case,
however the adaptation of the proof needs some changes. Select a smooth connected curve 
$\gamma_i$ in $\Delta$. We are going to perturb $F$ in a $\gamma$-neighborhood of $\Delta$.
By Lemma \ref{disjoint}, the perturbations can be made in each connected component
$\gamma_i$ in an independent way.

Recall that, for $k$ large enough, the curve $\gamma_i$ is contact, i.e. $\forall x\in 
\gamma_i, ~ T_xC=T_x \gamma_i \oplus D_x$. Moreover, the angle between $T_x \gamma_i$ and 
$D_x$ is bounded below by a uniform constant because of the transversality of the sequence.
Using the contact
metric $g_J$ associated to the fixed complex structure $J$ to define a geodesic flow, we 
can obtain a diffeomorphism:
$$ \phi_i: U_{\rho} \to V_{\rho} \subset S^1\times \C^n, $$
where $U_{\rho}$ is the $\rho$-neighborhood of $\gamma_i$ (in $g_k$ metric) and $V_{\rho}$
is its image by the flow, which is an open neighborhood of $S^1\times \{ 0 \}$. 
We can construct a metric in $S^1$ 
by imposing the condition that $(\phi_i)_{|\gamma_i}$ is an isometry with respect to the
rescaled contact metric $g_k$. In $\C^n$ we will fix the standard metric. 
The product metric will
be denoted by $g_0^k$. We can select with this choice a uniform $\rho>0$ such that
\begin{equation}
\lambda_m g_k(d\phi_i(v), d\phi_i (w)) \leq |g_0^k(v,w)| \leq \lambda_M g_k(d\phi_i(v),
d\phi_i(w)), ~~ \forall v,w \in T_x C, \forall x\in U_{\rho}, \label{casi_iso}
\end{equation}
where $\lambda_m,\lambda_M>0$ are uniform constants. Once we have fixed $\phi_i$, we obtain
a distribution $D_k$ in $\phi_i(U_{\rho})$ constructed as the image of
the distribution $D$. Denote by $D_h$ the integrable distribution given as $\{ p \}
\times \C^n$ defined in $S^1\times \C^n$. Perhaps after shrinking $\rho$ uniformly, we can 
check that
\begin{equation}
\angle_M(D_k(s,z), D_h(s,z))< c_u|z|k^{-1/2}, ~ \forall (s,z)\in V_{\rho}
\subset S^1\times \C^n, \label{apr_dist}
\end{equation}
where $c_u>0$ is a uniform constant. Moreover we can impose without loss 
of generality that $(\phi_i)_* J(\gamma_i)= J_0$. So we can project orthogonally
$(\phi_i)_*J$ (defined on $D_k$) to $D_h$ obtaining a new almost complex
structure $\hat{J}$. In fact it is easy to check, as in \cite{Do96}, that
\begin{equation}
\partial=\partial_0+\bar{\mu}\bar{\partial}_0, ~~ \bar{\partial}=\bar{\partial}_0+
\mu\partial_0, \label{apr_cplx}
\end{equation}
where $\partial$ and $\partial_0$ are the operators defined by the structures $\hat{J}$
and $J_0$ in $S^1\times \C^n$, and $|\mu(z)|\leq c|z|k^{-1/2}$, where $c$ is a uniform
constant. In order to finish the proof we follow these steps: we will define the perturbation, 
afterwards we will prove 
it satisfies the conditions for the distribution $D_h$ with the almost complex 
structure $\hat{J}$, and finally we will check the result for the distribution $D_k$.

Given any differentiable function $f:C \to \C$, we denote $f_0=f \circ \phi_i^{-1}$. By 
the inequalities (\ref{casi_iso})
we can use $F_0$ instead of $F$ for all the computations using the induced distribution $D_k$.
To construct the perturbation we define the complex Hessian $H=\frac{1}{2}\partial
\partial F$. Using the trivialization $\phi_i$ we may regard it as
$$ H_0(s,z)=\sum H_{\alpha \beta}(s) z_{\alpha}z_{\beta}, $$
on $\phi_i(U_{\rho})$. Also we take a cut-off function $\beta_{\rho}: S^1\times \C^n \to
[0,1]$ satisfying
\begin{enumerate}
\item $\beta_{\rho}(\phi_i(p))=1$, if $d_k(p,\gamma_i)\leq \frac{\rho}{2}$.
\item $\beta_{\rho}(\phi_i(p))=0$, if $d_k(p, \gamma_i)\geq \rho$.
\item $|\nabla \beta_{\rho}|=O(\rho^{-1})$.
\end{enumerate}
We can adjust $\beta$ to assure condition 3 because of equation (\ref{casi_iso}). 
The constant $\rho<\rho_0$ will be fixed along the proof to assure that the 
conditions are satisfied (namely we will have to shrink $\rho$ in a uniform way).
A modification of $F$ will be
$$ f_0(s,z)=\beta_{\rho}(w'(s)+H_0(s,z))+(1-\beta_{\rho})F_0(s,z), $$
where $w':S^1 \to \C$ is any smooth function. 

We denote by $\partial$ and $\bar{\partial}$ the operators associated to the almost complex 
structure $\hat{J}$ acting on $D_h$. Respectively we denote by $\partial_k$ and 
$\bar{\partial}_k$ the operators associated to the action of $(\phi_i)_*J$ in $D_k$.
The equivalent to Lemma 10 of \cite{Do99} is
\begin{lemma}
If $\rho>0$ is small enough, $k$ is sufficiently large, $|w'(s)-F_0(s,0)|$ is sufficiently
small and $|\frac{d w'(s)}{ds}|=O(1)$, then the inequality 
$|\partial_k f_0|\leq |\bar{\partial}_k f_0|$ is only satisfied in $\gamma_i$.
\end{lemma}
\noindent {\bf Proof:}
First, assume that we are at a point where $\beta_{\rho}=1$. Then $f_0=w'+H_0$ and 
$$ \partial f_0= \partial H_0, ~~ \bar{\partial} f_0=\bar{\partial} H_0. $$
The $\eta$-transversality of $\partial F_0$ yields the bound
$$ |\partial H_0(s,z)| \geq \eta |z| - |\bar{\partial}(\partial F_0)_{(z=0)}||z|. $$
Now recall that $\bar{\partial} \partial+ \partial \bar{\partial}=0$ on functions
and that the norm of $\partial \bar{\partial} F_0$ is controled because $F$ has uniform
mixed $C^3$-bounds. Then we can write
$$ |\partial H_0(s,z)| \geq \eta |z| - c_uk^{-1/2}|z|. $$
Using the inequality (\ref{apr_dist}) and that  $|\frac{\partial}{\partial s} 
H_{\alpha \beta}(s)|=O(1)$ (as follows from the mixed $C^3$-bounds of $F$) and 
recalling the bound $|\frac{dw'}{d s}|=O(1)$, we rewrite the inequality as
\begin{equation}
|\partial_k H_0(s,z)| \geq \eta |z| - c_u'k^{-1/2}|z|, \label{upper}
\end{equation}
where $c_u'>0$ is a uniform constant. On the other hand,
$$ |\bar{\partial} H_0| \leq c_u|z|^2k^{-1/2}, $$
and thus, by an analogous argument,
\begin{equation}
|\bar{\partial}_k H_0| \leq c_u''|z|^2k^{-1/2}+c_u''|z|k^{-1/2}. \label{lower}
\end{equation} 
Now, if we impose that $|\partial_k H_0|\leq |\bar{\partial}_k H_0|$, we obtain by comparing
(\ref{upper}) and (\ref{lower}), $z=0$ for $k$ large enough.

Now we study points in the annulus containing the support of $\nabla \beta_{\rho}$. In this
case,
$$\bar{\partial} f_0 =\bar{\partial} \beta_{\rho}(w'+H_0-F_0)+\beta_{\rho}\bar{\partial} 
H_0+(1-\beta_{\rho})\bar{\partial} F_0. $$
Bounding the right hand side as in \cite{Do99} we obtain an expression for the value 
$|\bar{\partial} f_0|$. Using
(\ref{apr_dist}) and the bounds $|\frac{d w'}{ds}|=O(1)$ and $|\bar{\partial}_k F_0|
=O(k^{-1/2})$, we conclude
$$ |\bar{\partial}_k f_0| \leq c(\rho^2+k^{-1/2}+|F_0(s,0)-w'|\rho^{-1}). $$
In the same way we know that
$$\partial_k f_0 =\partial_k \beta_{\rho}(w'+H_0-F_0)+\beta_{\rho}\partial_k H_0
+(1-\beta_{\rho})\partial_k F_0. $$
Using the transversality of $F$ we can obtain a lower bound for $|\partial_k f_0|$. The 
argument follows the one of Donaldson arriving to the final expression
$$ |\partial_k f_0|-|\bar{\partial}_k f_0| \geq \frac{\eta \rho}{2}-c(\rho^2+k^{-1/2}
+|w'-F_0(s,0)|\rho^{-1}). $$
Obviously, once fixed a sufficiently small $\rho$, for $k$ large enough and 
$|F_0(s,0)-w'|$ small enough compared with $\rho$ the inequality is strictly positive
for any point in the annulus.
\hfill $\Box$

To finish the proof we only have to check that the function $f$ conforms the local models
at any point of $\Delta$. We use the perturbation $w'$ to assure that the curves $\gamma_i$
project smoothly into $\CP^1$, this is equivalent to impose
\begin{equation}
\nabla_v f(x)\neq 0, \nonumber
\end{equation}
for any point $x\in \Delta$ and any nonzero $v\in T_x\Delta$. This can be achieved by a
generic perturbation, and so $w'$ can be selected to get it. Also we can assure, using 
this perturbation, that the intersections of branches of $f(\Delta)$ are transverse, again 
by a genericity argument. 

Finally, given any $x\in \Delta$ there exists coordinates $(s,x_1, y_1, \ldots, x_n, y_n)$
such that $f$ is locally written as
$$ f(s,z_1, \ldots, z_n)= \varphi(s)+\sum H_{\alpha \beta}(s)z_{\alpha} z_{\beta}, $$
where $\varphi'(0)\neq 0$ (this condition is equivalent to the smoothness of the branches
of $f(\Delta)$).
Now at a neighborhood of $s=0$ the $1$-parametric family $(H_{\alpha \beta}(s))$ of
bilinear complex forms can be diagonalized by a smooth family of complex changes of 
coordinates if
the eigenvalues of $(H_{\alpha \beta})$ are all distint. This is a genericity condition
that can be achieved at all the points of $\Delta$ by a generic perturbation of $O(k^{-1/2})$
in $F$, before starting the perturbation process which we have developed along this section.
With this condition, we obtain a smooth family of invertible complex matrix $P(s)$ such that
$$ \sum H_{\alpha \beta}(s)z_{\alpha} z_{\beta}= z^TP(s)^TP(s)z, $$
where $z=(z_1, \ldots, z_n)$. Therefore the change of coordinates
\begin{eqnarray*}
\R\times \C^n & \to & \R\times \C^n \\
(s,z) & \to & (s, P(s)z)
\end{eqnarray*}
gives us the required local model. This finishes the proof of the main theorem.

\section{Transversality results.} \label{asymp}
In this section we will prove Theorem \ref{trans_2}. 
\subsection{The globalization scheme}
First we recall how the globalization process developed in \cite{Do96} adapts to 
the contact setting. This adaptation was carried out in \cite{IMP99}. Now, we set up the 
process in a functorial way in the style of \cite{Au99}. As always, we denote by $(C,\theta)$
an exact contact manifold.
\begin{definition} \label{open_prop}
A family of properties $\CrP(\epsilon,x)_{x\in C, \epsilon>0}$ of sections of bundles
over $C$ is local and mixed $C^r$-open if, given a section $s$ satisfying $\CrP(\epsilon,x)$,
any section $\sigma$ such that $s-\sigma$ has $(\eta, c_R)$ mixed $C^r$-bounds satisfies
$\CrP(\epsilon-c_u\eta, x)$,  for some constant $c_u$.
\end{definition}
\begin{proposition}[\cite{IMP99}] \label{global}
Let $\CrP(\epsilon,x)_{x\in C, \epsilon>0}$ be a local and mixed $C^r$-open family of
properties of sections of vector bundles $E_k$ over $C$. Assume there exist uniform
constants $c$, $c'$, $c''$, $p$ and a function $f:\R^3 \to \R^+$ such that, given any $x\in 
C$, any small enough $\delta>0$, and mixed $C^r$-bounded sections $s_k$
of $E_k$ with uniform mixed $C^r$-bounds, say $(c_D,c_R)$, there exist, for all large 
enough $k$, mixed $C^r$-bounded sections $\tau_{k,x}$ of $E_k$ with the 
following properties:
\begin{enumerate}
\item $\tau_{k,x}$ has mixed $C^r$-bounds $(c''\delta, f(c''\delta,c_D,c_R))$,
\item the sections $\frac{1}{\delta}\tau_{k,x}$ have mixed Gaussian decay away from
$x$ in $C^r$-norm,
\item $s_k+\tau_{k,x}$ satisfy the property $\CrP(\eta,y)$ for all $y\in B_{g_k}(x,c)$,
with $\eta=c'\delta(\log (\delta^{-1}))^{-p}$.
\end{enumerate}

Then, given any $\alpha>0$ and mixed $C^r$-bounded sections $s_k$ of $E_k$, 
there exist, for all large enough $k$, mixed $C^r$-bounded sections $\sigma_k$
of $E_k$, such that $s_k-\sigma_k$ has mixed $C^r$-bounds $(\alpha, c_R)$ for some
$c_R>0$. Also, the sections $\sigma_k$ satisfy $\CrP(\epsilon,x)$ for some uniform 
$\epsilon>0$ at any $x\in C$.
\end{proposition}
{\bf Sketch of the proof:}
Although this is just a slight variation of the globalization argument in \cite{IMP99}
we provide the main lines by completeness.
Let $S$ be a finite set of points in $C$ verifying the following properties:
\begin{enumerate}
\item $\bigcup_{x\in S} B_{g_k}(x,c)\supset C$.
\item There exists a partition $S=\bigcup_{j \in J} S_j$ verifying that $d_{g_k}(x,y)>\CrN$
if $x,y\in S_j$. $\CrN$ will be fixed along the proof.
\item The cardinal of $J$ is $O(\CrN^{2n+1})$.
\end{enumerate}
The idea is to achieve the property $\CrP(\epsilon,x)$ in all the balls $B_{g_k}(x,c)$, 
$x\in S_j$  at once. Using the hypothesis at each $x\in S_j$ we can build a local section 
$\tau_{k,x}$ which satisfies $\CrP(\eta,y)$ for all $y\in B_{g_k}(x,c)$, where $\eta=
c''\delta(\log (\delta^{-1}))^{-p}$. We select $\CrN$ large enough, such 
that at a given $x\in S_j$ the components of the perturbations due to the rest of the points of
the set $S_j$ through the directions defined by the distribution $D$ do not destroy the 
transversality obtained at $B_{g_k}(x,c)$. The condition we impose to $\CrN$ is
\begin{equation}
c''\delta(\log (\delta^{-1}))^{-p} \geq 2c_uc'\delta \exp(-\lambda\CrN^2), \label{key_in}
\end{equation}
where $c_u$ is a uniform constant. The left hand side of the inequality is the amount of 
transversality obtained by the perturbation at a point $x\in S_j$. The right hand side is a 
value less than the double of
the value of the norm of the sum of all the other perturbations in the set $S_j$ (different 
from the selected and multiplied by the uniform constant provided by the local
an $C^r$-open property). So the inequality assures us that, after adding the perturbation,
we have obtained $\CrP(\frac{1}{2}c''\delta(\log (\delta^{-1}))^{-p},y)$ at a 
$c$-neighborhood of $S_j$. 

Now the process is clear. We perturb in $S_1$ to obtain $\CrP(\eta_1,x)$
at a $c$-neighborhood of it. Afterwards we add again a perturbation with mixed $C^r$-bounds 
$(\frac{1}{c_u}\eta_1/2, c_R^2)$ to 
achieve $\CrP(min(\eta_2,\eta_1/2),y)$ for all $y$ at a $c$-neighborhood of $S_1\bigcup S_2$. 
The number of steps is independent of $k$, so the achieved final property is $\CrP(\eta,x)$
where $\eta$ does not depend on $k$. The added section has also mixed $C^r$-bounds
$(\epsilon,c_R)$. The constant $\epsilon$ is obtained by adjusting the $\eta_i$ in the 
iteration, which is possible by property 1. The constant $c_R$ is uniform because at each 
stage does depend only on the precedent values of $c_R$ and $\eta_j$, so it is independent 
of $k$, and thus the uniformity is obvious since the number of 
stages is independent of $k$ (i.e. $O(\CrN^{2n+1})$). Remark that $c_R$ cannot be bounded
since property 1 provides uniformity but not control of the constant!

To end we have to check the inequality
(\ref{key_in}) in all the steps of the process. Following the asymptotic analysis of 
\cite{Do96} we conclude that it is possible to obtain an integer $\CrN$ independently of the 
step $j$. This ends the proof. 
\hfill $\Box$

\subsection{The local perturbation.}
We are going to achieve the transversality property for $\partial F$, required in
Theorem \ref{trans_2}, by using Proposition
\ref{global}. The first observation is the following
\begin{lemma} \label{lejos}
There exists a uniform constant $\xi>0$, such that, for $k$ large enough, any point $x$
verifying $|\partial F(x)|<\xi$ lies in the set $\Omega_{2\xi}$.
\end{lemma}
The proof is identical to that of Lemma \ref{far} and we refer again to \cite{Do99} for
details. 

We are going to perturb the sequence of sections $s_k^1$ on order to obtain
transversality.
The open mixed $C^2$-open property  $\CrP(\epsilon,x)$ which we have to 
obtain in $C$ is that $\partial F$ is 
$\epsilon$-transverse to {\bf 0} at $x$. We are in the hypothesis of 
Definition \ref{open_prop}. Now we are going to construct a local section satisfying the 
hypothesis of Proposition \ref{global}, thus concluding the proof. 
We can choose $c>0$ uniformly to assure that $B_{g_k}(x,8c)\subset \Omega_{\xi}$ for
any $x\in \Omega_{\frac{3}{2}\xi}$ and also that $B_{g_k}(x,8c)$ is in the complementary
set of $\Omega_{\xi}$ for any $x$ in the complementary set of $\Omega_{\frac{3}{2}\xi}$.

If we are in the complementary set
of $\Omega_{\frac{3}{2}\xi}$ we can choose the section $\tau_{k,x}$ to be zero, because of 
Lemma \ref{lejos}.

Choose a point $x\in \Omega_{\frac{3}{2}\xi}$. Take an approximately contact holomorphic chart
$(s, z_k^1, \ldots, z_k^n)$ satisfying the properties of Corollary \ref{approx}. Generalizing 
\cite{Au99,MPS99} we define the $1$-forms
$$ \mu_k^j=\partial(\frac{z_k^j s_{k,x}^{ref}}{s_k^0}). $$
In $\Omega_{\xi}$ these forms have uniform mixed $C^3$-bounds because the term 
$\nabla^r \bar{\partial}
z_k^j$ can be bounded near $x$, and furthermore both $s_{k,x}^{ref}$ and $s_k^0$ have 
mixed $C^3$-bounds, as well as $s_k^0$ is also bounded below. Recall that 
$(\mu_k^j)_{j=1,\ldots, n}$
is a unitary basis of $D^*$ at $x$, and it is almost unitary in the ball $B_{g_k}
(x,8c)$ (i.e. the basis $\mu_k^j$ is arbitrarily close to a unitary basis), after 
eventually shrinking $c$ uniformly. We define an application
$v:B_{g_k}(x,8c) \to \C^n$ by the formula:
$$ \partial F= \sum_{j=1}^n v_j\mu_k^j. $$
In matrix notation,
\begin{equation}
\partial F= v^T\cdot \mu_k, \label{helping}
\end{equation}
where $\mu_k^T= (\mu_k^1, \ldots, \mu_k^n)$. Thus it is possible to understand $\mu_k $ as
a linear map $\mu_k(y): \C^n \to T^*_yM$. Multiplying by $\mu_k^{-1}$ in (\ref{helping})
we obtain
$$ v^T= \partial F \cdot \mu_k^{-1}. $$
This implies that $|v|=O(1)$, since $\mu_k$ is approximately unitary. 
We compute now the derivatives of $v$ using equation (\ref{helping}). To do this,
recall that $\partial F$ has mixed $C^2$-bounds in the ball $B_{g_k}(x,8c)$ (because
$F$ has mixed $C^3$-bounds). Also, in the same way, $\mu_k$ has mixed $C^2$-bounds.
Differentiating (\ref{helping})
$$ \nabla \partial F= \nabla v^T \cdot \mu_k + v^T \cdot \nabla \mu_k, $$
which implies that $|\nabla v^T \cdot \mu_k|=O(1)$ and using again that $\mu_k$ is
approximately unitary, we finally get
$$ |\nabla v|= O(1) $$
Differentiating respect to $\bar{\partial}$ we find that $|\bar{\partial}
v|=O(k^{-1/2})$, and iterating the process
$$ |\nabla \nabla v|=O(1), ~~ |\nabla \bar{\partial} v|= O(k^{-1/2}). $$
Now we use the approximately contact-holomorphic chart $\Psi$ defined in
Corollary \ref{approx}, after eventual uniform shrinking of $c$. We construct the function 
$\hat{v}=v \circ \Psi^{-1}$. By the properties of $\Psi$, it is easy to check that:
$$ |\hat{v}|=O(1), ~~ |\nabla^r \hat{v}|=O(1), ~~ | \nabla^{r-1} \bar{\partial}
\hat{v}|=O(k^{-1/2}), ~~ r=1,2. $$
Scaling the coordinates by a uniform constant we can assume that $\Psi(B_{g_k}(x,2c))
\subset B_{2n+1}(0,2) \subset \Psi(B_{g_k}(x,8c))$. Now, we are in the hypothesis
of the following 
\begin{proposition} \label{trans_cont}
Let $f_k\colon B \times [0,1] \to \C^m$ be a sequence of functions where $B$ 
is the ball of radius $1$ in $\C^n$ and $B\times [0,1]$ is equipped with 
a sequence of contact forms $\theta(k)$ whose distributions are asymptotically flat. 
Let $0 < \delta < 1/2$ 
be a constant and let $\sigma = \delta (\log (\delta^{-1}))^{-p}$, where $p$ is a
integer depending only on the dimensions. Assume that $f_k$ satisfies over $B\times  [0,1]$
the following bounds
$$ |f_k| \leq 1, ~~~~ |\bar{\partial}_0 f_k| \leq \sigma, ~~~~|\nabla
\bar{\partial}_0 f_k| \leq \sigma,$$
for $k$ large enough, where $\bar{\partial}_0$ is the $(0,1)$ operator defined in $D(k) = \ker
\theta(k)$ by vertical projection of the standard complex structure $J_0$.  Then 
for $k$ large enough
there exists a smooth curve $w_k\colon [0,1] \to \C^m$  such that $|w_k|< \delta$
and the function $f_k - w_k$ is $\sigma$-transverse to zero on $B(0,1/2)\times [0,1]$. 
Moreover, if $|\partial f_k /\partial s| < 1$ and $|\partial \nabla f_k /\partial
s| < 1$, we can choose $w_k$ such that $|d^iw_k /ds^i | < \Phi(\delta)$, $(i=1,2)$; $d^jw_k /ds^j (0) = 0$ and 
$d^jw_k / ds^j (1) = 0$, for all $j \in \N$, where $c$ is a uniform  constant and $\Phi:\R^+
\to \R^+$ is a function depending only on the dimensions.
\end{proposition}
This proposition will be proved in Section \ref{local_res}. This is
the analogous in the contact case to Theorem 12 in \cite{Do99}. For the particular
value of $m=1$ it has been proved in \cite{IMP99}. 

Now we apply this proposition to the map $\hat{v}$ over $[-1,1]\times B$, 
because for $k$ large enough,
it satisfies the hypothesis (without loss of generality we can choose $[-1,1]\times
B$ instead of $[0,1]\times B$, also we suppose $|\hat{v}|_{C^2}\leq 1$, multiplying by
uniform constants). The obtained path $\hat{w}$ is extended to $\R \times \C^n$
as
\begin{equation}
\hat{w}(s,z)= \left\{ \begin{array}{ll} \hat{w}(1) & , \hbox{\rm for } s>1, \\ \hat{w}(s) & 
, \hbox{\rm for } s\in [-1,1], \\ \hat{w}(-1) & , \hbox{\rm for } s<-1. \end{array} \right. 
\nonumber
\end{equation}
We keep the same notation for this extended map. The map $\hat{v}-\hat{w}$ is transverse 
to {\bf 0} in $B(0,1)\times [-1,1]$, thus $\hat{v}-\hat{w}$ 
is also transverse to {\bf 0} in $B_{2n+1}(0,1)$. Remark that the chart $\Psi$ is defined
in a ball of radius $O(k^{1/2})$ . So the pull-back of $\hat{w}$ by $\Psi$, denoted $w$, is
well defined in this ball, and it will be enough for our purposes. By the properties of $\Psi$
we have that $v+w$ is $\eta$-transverse to {\bf 0} over the ball $B_{g_k}(x,c)$. The 
constant $\eta$ does not coincide with the transversality obtained by $\hat{v}+
\hat{w}$ by a uniform factor, so $\eta=c_u\delta (\log (\delta^{-1}))^{-p}$, the constants
$c_u$ and $p$ are uniform and $\delta$ is the norm of $w$ except by a uniform factor
$c_{u'}$, i.e. $|w|<c_{u'}\delta$, $c_{u'}\neq 0$. 

The needed perturbation is
$$ \tau_{k,x}= \sum_{j=1}^n w_j z_j s_{k,x}^{ref}. $$ 
The first question is whether $\tau_{k,x}$ has the desired mixed $C^3$-bounds. The 
verification is mere routine and it follows the lines of \cite{IMP99}. As an example, we 
compute the bound for $\bar{\partial} \tau_{k,x}$:
$$ \bar{\partial} \tau_{k,x}= \sum \bar{\partial} w_j z_j s_{k,x}^{ref}+
w_j \bar{\partial} z_j s_{k,x}^{ref} + w_j z_j \bar{\partial} s_{k,x}^{ref}. $$
The third term is easily bounded by
\begin{eqnarray*}
 |w_jz_j \bar{\partial} s_{k,x}^{ref}|& \leq & O(1)2d_k(x,y) c_ak^{-1/2}P(d_k(x,y))
\exp (-\lambda d_k(x,y)^2)= \\
& = & k^{-1/2}Q(d_k(x,y))\exp(-\lambda d_k(x,y)^2).
\end{eqnarray*}
The second one is bounded in the same way recalling that $\bar{\partial} z= c_0k^{-1/2}
d_k(x,y)$, where $c_0>0$ is uniform. For the first term we proceed as follows:
$$ |\bar{\partial} w_j z_j s_{k,x}^{ref}| \leq |\bar{\partial} w_j| P(d_k(x,y))
exp(-\lambda d_k(x,y)^2). $$
Recall that $\theta_{k^{1/2}}$ is asymptotically flat. Using this fact we obtain that
$|\nabla w_j|\leq O(k^{-1/2})d_k(x,y)$, thus in particular $|\bar{\partial} w_j| 
\leq O(k^{-1/2}) d_k(x,y)$. This concludes the bounding.

Now we have to compare $\tau_{k,x}$ with $w$. If
$\partial \frac{\tau_{k,x}}{s_k^0}$ were $w^T\mu_k$, the proof would be rapidly
concluded. But we are almost in this situation because
$$ \partial \frac{\tau_{k,x}}{s_k^0}=
\Sigma \partial w_j\frac{z_js_{k,x}^{ref}}{s_k^0}+ \Sigma w_j\mu_j. $$
Using (\ref{apr_dist}) again, we obtain
$$ \partial \frac{\tau_{k,x}}{s_k^0}= O(k^{-1/2})+w^T \mu_k. $$
Therefore, we can find $\tilde{w}$ such that $\partial \frac{\tau_{k,x}}{s_k^0}=\tilde{w}^T
\mu_k$ and $|w-\tilde{w}|=0(k^{-1/2})$. So, for $k$ large enough $v-\tilde{w}$
is, say,  $0.9\eta$-transverse to {\bf 0} in $B_{g_k}(x,c)$.
We end by remarking that $\partial \frac{s_k^1-\tau_{k,x}}{s_k^0}$ is $c_b\eta$-transverse
to {\bf 0}, $c_b>0$ a uniform constant, if $v-\tilde{w}$ is $0.9\eta$-transverse to
{\bf 0}. This is obvious since both just differ by the application of the almost unitary 
matrix $\mu_k$.
\hfill $\Box$

\section{Local results.} \label{local_res}
The aim of this section is to prove Theorem \ref{trans_cont}. This is the 
generalization of the local results in the symplectic setting, needed to achieve
controled transversality in the contact case.
\subsection{Reduction to integral distributions.}
We can easily reduce the proof to the following
\begin{proposition} \label{par_trans_cont}
Let $f\colon B \times [0,1] \to \C^m$ be a complex valued function, where $B$ 
is the ball of radius $1$ in $\C^n$. Let $0 < \delta < 1/2$ be a constant 
and let $\sigma = \delta (\log (\delta^{-1}))^{-p}$, where $p$ is a suitable fixed
integer depending only on the dimensions $n,m$. Assume that $f_s$ satisfies the following 
bounds over $B\times [0,1]$
$$ |f_s| \leq 1, ~~~~ |\bar{\partial}f_s| \leq \sigma, ~~~~|\nabla 
\bar{\partial}f_s| \leq \sigma.$$
Then there exists a smooth curve $w \colon [0,1] \to \C^m$ 
such that $|w|< \delta$ and
the function $f_s - w(s)$ is $\sigma$-transverse to zero over the  ball $B(0,1/2)$. 
Moreover, if $ |\partial f_s/\partial s| < 1$ and  $|\partial
\nabla  f_s /\partial s| < 1$, we can choose $w$ such that
$|d^iw /ds^i| < \Phi(\delta)$ ($i=1,2$),
$d^jw /ds^j(0) = 0$ and $d^jw /ds^j(1) = 0$ for all $j\in \N$,
where $\Phi$ is a function depending only on the dimensions $n,m$.
\end{proposition}
The proof of this Proposition is a generalization of Lemma
10 in \cite{IMP99} to the case $m>1$. 

{\bf Proof of Theorem \ref{trans_cont}:}
We have to obtain the transversality of $f_k-w_k$
when we restrict $\nabla (f_k-w_k)$ to the distribution defined by $\theta(k)$. 
Recall that $\theta(k)$ is asymptotically flat, so $\angle_M(D_h, D_k)=O(k^{-1/2})$,
where $D_h=\C^n \times \{p \}$ and $D_k=\Ker ~ \theta(k)$.
The key idea is that $\nabla (f_k-w_k)=O(1)$, so if $f_k-w_k$ is $\eta$-transverse to
$D_h$, then it is, say, $0.9\eta$-transverse to $D_k$ for $k$ large enough. The factor
$0.9$ can be eliminated by increasing $p$ uniformly.
\hfill $\Box$

\subsection{Proof of Proposition \ref{par_trans_cont}.}
Let $f:\C^n \to \C^m$ be a smooth application. We define the subset $U(f,w,\delta, \sigma)$ of
the ball $B(w,\delta)$ of radius $\delta$ as the set of points such that $f$ is 
$\sigma$-transverse to $w$ in $\frac{1}{2}B^{2n}$, the ball of radius $1/2$ in $\C^n$.
First we prove the following
\begin{theorem}[Extension of Theorem 12 in \cite{Do99}] \label{main_key}
For any $n$, $m$, $\delta>0$ and $0<\gamma<1$ there is a $p=p(n,m,\gamma)$ such that if
we define $\sigma=\delta (\log (\delta^{-1}))^{-p}$ then for 
all the maps $f: B^{2n} \to \C^m$ verifying that
$$ |f|\leq 1,~~ |\bar{\partial} f|\leq \sigma, ~~ |\nabla \bar{\partial} f|\leq \sigma, $$
there exists a connected component of $U(f,w,\delta, \sigma)$ containing
another path-connected set $U'(f,w,\delta, \sigma)$ whose volume is at least
$\gamma$ times the total volume of $B(w,\delta)$ and such that given two points
$x,y\in U'(f,w,\delta,\sigma)$ is possible to find a smooth curve $\gamma$ in
$U(f,w,\delta, \sigma)$ joining $x$ and $y$ with curvature at each point and length
bounded by $\Phi(\delta)$, where $\Phi$ is a function depending only on the dimensions $n,m$.
\end{theorem}

We repeat the the proof of Theorem 12 of \cite{Do99} taking care of 
some details. 

\noindent {\bf Proof:}
The first step in the proof is to approximate $f_s$ by a holomorphic function $\hat{f}_s$ such
that $|f_s-\hat{f}_s|_{B,C^1}<c\sigma$ (see Lemma 28 in \cite{Do96}). This process 
does not hold over all of the unit ball. This is the reason why we restrict
ourselves to the ball $B'=\frac{1}{2}B^{2n}$ in the sequel. Then we approximate $\hat{f}_s$
by a polynomial. We can obtain polynomials $g_s$ such that $|g_s-f_s|_{B', C^1}\leq c\sigma$
and their degree $d$ can be stimated by $O(\log(\sigma^{-1}))$.

Adapting notations of \cite{Au97, IMP99} we denote by $Z_{h_s,\epsilon}$ the 
images of the set of points of $B'$ which are not $\epsilon$-transverse to
{\bf 0} for $h_s$. We want to prove that one component of the complementary set of $Z_{f_s,
\sigma}$ satisfies the required properties. The first observation is that the $C^1$-closedness
of $f_s$ and $g_s$ assures us that $Z_{f_s, \sigma} \subset Z_{g_s, (c+1) \sigma}$.

We use now the following
\begin{theorem}[Theorem 26 of \cite{Do99}] \label{Donald_help}
Given a polynomial $g:\C^n \to \C^m$ of degree $d$ and $\epsilon>0$ there is a real-algebraic
subvariety $A(g)\subset \C^m$ of codimension $2$ and degree $D$ such that $Z_{g,\epsilon}$
is contained in the $K\epsilon$-neighborhood of $A(g)$, where $K,D\leq (d+1)^p$, for some
integer $p$ depending only on the dimensions $n,m$.
\end{theorem}

Also we use the following Donaldson's result:
\begin{proposition} [Proposition 31 in \cite{Do99}] \label{peque}
For each integer $N$ and real number $\theta>0$, there is a $\mu=\mu(\theta,N)$
with the following property. For any real-algebraic hypersurface $A\subset \R^n$ of
degree $D$ and $\epsilon\leq (D+1)^{\mu}$,
$$ Vol(B^N\cap A_{\epsilon}) \leq \theta.$$
\end{proposition}

Denote $\sigma'=(c+1)\sigma$. With these two results and following the discussion in 
\cite{Do96} p. 689 about the behaviour
of the function $\delta (\log \delta^{-1})^{-p}$ we can assure the 
$\sigma'=\delta \log (\delta^{-1})^{-p}$-neighborhood, $p$ a fixed integer, of the bad set 
$A(g_s)$ has volume arbitrarily small. Also we can assure the same condition for the 
$3\sigma'$-neighborhood (changing $p$ slightly). We take a covering of $B(0,\delta)$ by
balls $B(x_i,\sigma'/2)$ of centers $x_i$ and radius $\sigma'/2$ and assuring that the 
covering of balls with radius $\sigma'$ centered in the same points cover each point of 
$B(0,\delta)$ only a finite uniform number of times, for instance less than $\nu$ times. 
Denote by $\CrC$ the
set of centers $x_i$ of the balls of the covering contained in the $2\sigma'$-neighborhood 
of $A(g_s)$. Recall that the union of these balls is contained in $A_{3\sigma'}(g_s)$,
so we can conclude that 
$$ \sum_{x_i\in Cr_C} vol(B(x_i,\sigma')) \leq \nu vol(A_{3\sigma'}(g_s)).$$
And from this expression we can easily obtain a bound for the number of balls in $\CrC$ as
$$ card(\CrC)\leq c_u \delta/ (\sigma')^2, $$
where $c_u$ is a uniform constant and so this number only depends on $\delta$. But observe
that 
$$ A_{\sigma'}(g_s) \subset \bigcup_{x_i \in \CrC} B(x_i, \sigma'/2) \subset 
\bigcup_{x_i \in \CrC} B(x_i, \sigma') =W_s \subset A_{3\sigma'}(g_s) $$.

Recalling that $W_s$ has an arbitrary small volume, a standard isoperimetric inequality
gives us that in the complementary of $W_s$ we can find a connected component of
arbitrarily big volume (see, for instance, \cite{Au99} for more details). This
component will be the set $U'(f,w,\delta, \sigma)=U'$ in the statement of the Theorem.
Obviously in the complementary of $A_{\sigma}(f_s)\subset A_{\sigma'}(g_s)$
there will be a connected component containing $U'$, this component will be the set
$U(f,w,\delta, \sigma)=U$. To finish we have only to check that these two sets
satisfy the required properties. 

We adapt the ideas of \cite{IMP99}. We call $N=card(\CrC)$.
Observe that we have fixed $\sigma$, and then $N=f(\delta)$ and $2\pi\sigma N=g(\delta)$, 
for some functions $f$ and $g$. Now take $y,z$ points in the large connected component $U$ 
of the complementary of $Z_s^+$. Denote $L(y,z)$ the straight segment joining them. This 
segment cuts at most at $2N$ points $y_0,z_0,y_1,\ldots$ to the border 
$Bor=\partial(\bigcup_{x_i \in \CrC} B(x_i, \sigma'))$.
Obviously $L(z_i,y_{i+1})\subset U$ and $L(y_i,z_i)
\subset \bigcup_{x_i \in \CrC} B(x_i, \sigma')$. We replace the lines $L(y_i,z_i)$ 
by curves $C(y_i,z_i)$
contained in $Bor$ connecting $y_i$ and $z_i$. We construct the curves 
following maximal diameters of the spheres which define the border and
so $\length(C(y_i, z_i)) \leq g(\delta)$. Therefore the curve
$$\gamma'=L(y,y_0) \bigcup C(y_0,z_0) \bigcup L(z_0,y_1) \cdots$$ 
satisfies
\begin{eqnarray*}
\length(\gamma') & = & L(y,y_0) +L(z_q,z) + \Sigma L(y_i,y_{i+1}) + \Sigma C(y_i,z_i)   
\leq 2\delta + \\ 
& & + f(\delta)g(\delta) = \Phi(\delta)/2, 
\end{eqnarray*}
where $\Phi$ is some function depending only on the dimensions. Perturbing slightly
$\gamma'$ to make it differentiable and removing it from the border
we obtain $\gamma$ which, bounding enough the perturbation, satisfies $\length(\gamma)\leq
\Phi(\delta)$. So the length between two points can be bounded by a function of $\delta$.
Moreover, if we translate the diameters in the border $Bor$ till the border of
$\bigcup_{x_i \in \CrC} B(x_i, \sigma'/2)$, we can assure that the curvature of
the path can be bounded by $O(1/\sigma')$, again a function of $\delta$.

\hfill $\Box$ 

With the precedent result we can easily finish the proof of Proposition 
\ref{par_trans_cont}. Recall that $f_s$ satisfies
$$ |\frac{\partial f_s}{\partial s}|\leq 1, ~~ |\frac{\partial \nabla f_s}{\partial s}|\leq 1. 
$$
This implies by a symple application of the Mean Value Theorem of the Differential Calculus
(for vectorial applications) that:
$$ |f_s-f_{s+\epsilon}|< \epsilon, ~~ |\nabla f_s -\nabla f_{s+\epsilon}|<\epsilon. $$
Recall that given two applications $f,g:\C^n \to \C^m$ such that $|f-g|_{C^1}<\epsilon$,
if $f$ has a left inverse with norm less than $\eta^{-1}$, for a large enough uniform 
constant $c_0$, $g$ has a left inverse of norm less than $(\eta-c_0\epsilon)^{-1}$.
This implies that
\begin{equation}
U(f_s, 0, \delta, \sigma) \subset U(f_{s+\epsilon}, 0, \delta, \sigma-c_0\epsilon). 
\label{contr}
\end{equation}
We choose a number $\sigma$, which satisfies the hypothesis of Theorem
\ref{main_key}, with $\gamma>\frac{1}{2}$. Take an integer $q$ satisfying that
\begin{equation} 
\frac{1}{q}<\sigma/c_0. \label{paso}
\end{equation}
We can choose $x_i\in U(g_{i/q},0,\delta, (c+1)\sigma)\bigcap 
U(g_{(i+1)/q}, 0, \delta, (c+1)\sigma)$, for $i=0,\ldots, q-1$. Consider the path $H_i=x_i 
\times [i/q,(i+1)/q]$. There exists a smooth path $V_i$ connecting $x_i$ with
$x_{i+1}$ in $U'(g_{(i+1)/q}, 0, \delta, \sigma)$. Its length is bounded above 
by $\Phi(\delta)$, as well as its curvature at any point. We parametrize $V_i$ by its 
arch-length.

We choose a smooth function $\beta \colon [0,1] \to [0,1]$ satisfying:
$$
\beta(x)=  \left\{ \begin{array}{ll} 0, & x\in [0, 1/4], \\
0<\beta(x)<1, & 1/4<x<3/4, \\
1, & x\in [3/4, 1]. \end{array} \right.
$$
and compute $|\beta|_{C_2}=c_b$. Denote $\beta_i(x)=\beta(x)\cdot
\mbox{length}(V_i)$. This function has norm $|\beta_i|_{C_2}<c_b\Phi(\delta)$.
Define a path
$$ w((i/q+\epsilon) =
V_i\left( \beta_i \left( \frac{1}{q}\left( \epsilon + \frac{q}{2}\right)
\right) \right) , ~~~~  |\epsilon|\leq q/2. $$ 
The map $w$ is smooth and we obtain $|w|_{C_2}<c_b \Phi(\delta)$, a bound 
which is a function of $\delta$. Remark that the first derivative depends on the
length of the path $V_i$ and the second on its curvature. Finally we observe that 
$|V_i(s)-w(s)| < \frac{1}{2q}$. It implies, by
(\ref{contr}) and (\ref{paso}), that $w(s)$ is $\sigma/2$ transverse to {\bf 0}, for 
all $s\in[0,1]$. We can find an integer $p'$ such that 
$\sigma'=\delta(\log (\delta^{-1}))^{-p'}<\sigma/2$ and the proof is finished.

\section{Topological considerations.} \label{topology}
In this section we will do some simple topological remarks about the contact Lefschetz pencils.
In the first subsection we will study the topological relationship between the smooth
fibres of a contact fibration.
\subsection{Crossing critical curves.} \label{previous}
Over each of the connected components of $\CP^1-f_k(\Delta)$ the fibres of $f_k$ are isotopic.
Now, we are interested in analyzing the topological behaviour of the fibres when crossing 
through $\Delta$. We will prove in this subsection
the following
\begin{proposition} \label{handle}
Given a contact fibration $(f,\Delta, A)$, choose a path $\gamma:[0,1]\to \CP^1$ such that 
$\gamma(1/2)\in f(\Delta)$ and the rest of the points of $\gamma$ are regular values of $f$.
Then $N=f^{-1}(\gamma(1))$ is built up from $N'=f^{-1}(\gamma(0))$ by adding a $n$-dimensional
handle and removing another $n$-dimensional handle. Therefore $H_j(N)= H_j(N')$ (resp.
$\pi_j(N)=\pi_j(N')$) for $j=0, \ldots, n-2$.
\end{proposition}
{\bf Proof:}
We restrict ourselves without loss of generality to neighborhoods of $\gamma(1/2)$ and
of the critical point of $f$, where we 
can define a compatible chart. So, with the usual identifications, we can suppose 
that $f: \R \times \C^n \to \C$. Moreover, for simplicity, we will assume that 
$f(s,z_1, \ldots, z_n)= s+z_1^2+ \cdots +z_n^2$ (being the general case a 
straightforward generalization) and the path will be $\gamma(t)=2(t-1/2)i$ 
with a critical value for $f$ at $t=1/2$. 

The proof of the result reduces to show that $B=f^{-1}(\gamma([0,1]))$ is a cobordism between
$N=f^{-1}(\gamma(0))$ and $N'=f^{-1}(\gamma(1))$ with only one surgery of index $n$.
And this follows if we find a Morse function with a critical point of index $n$.
Choose $h= im(f)= 2x_1y_1+2x_2y_2+\cdots + 2x_ny_n$. We can assume at a neighborhood
$U$ of the critical point ${\bf 0}=(0, \ldots, 0)$ that $B\bigcap U \subset 
g^{-1}(0)$, where $g(s,x_1, \ldots)=Re(f)= s+x_1^2+y_1^2 +\ldots + x_n^2 + y_n^2$. To compute the index 
of $g$ we have to restrict ourselves to $\ker \nabla g({\bf 0})$ which is
$$\ker \nabla g({\bf 0})= \{ (0, z_1, \ldots, z_n): (z_1, \ldots, z_n)\in \C^n. \}$$
Finally, recall that
$$ \nabla \nabla h= \left( \begin{array}{ccc} 0 & 0 & 0 \\ 0 & 0 & 2I_n \\
0 & 2I_n & 0 \end{array} \right). $$
Restricting to $\ker \nabla g({\bf 0})$ we obtain
$$ \nabla \nabla h_{|\ker \nabla g({\bf 0})}=  \left( \begin{array}{cc} 0 & 2I_n \\
2I_n & 0 \end{array} \right), $$
which has index $n$.
\hfill $\Box$

Now we can obtain a geometrical relationship between the contact submanifolds obtained
in Theorem \ref{trans_1} as zero sets of transverse sequences with mixed $C^3$-bounds.
The result will be
\begin{proposition} \label{quasi_iso}
Let $S$ be a line bundle with connection over a closed exact contact manifold $(C,D)$
Given $N_k$ and $N_k'$ sequences of contact submanifolds obtained as zero sets of
transverse sections with mixed $C^3$-bounds of the bundles $S\otimes L^{\otimes k}$, 
then $N_k$ and $N_k'$ are cobordant through a cobordism defined by surgeries of index $n$.
\end{proposition}
{\bf Proof:}
Initially we will suppose that the contact form $\theta$ and the almost complex structure
$J$ used in the definition and construction of the sequences of sections are
coincident in the two cases.

Choose two of such sequences $\sigma_k^1$ and $\sigma_k^2$ which
are $\eta$-transverse to {\bf 0}, for some $\eta>0$. 
We claim that we can assume without loss of generality that $\sigma_k^1
\oplus \sigma_k^2$ is $\eta'$-transverse to {\bf 0}, for some $\eta>0$. If we suppose 
this to be false, we can apply Theorem \ref{trans_1} and find a sequence of sections
$\tau_k^1\oplus \tau_k^2$ with mixed $C^3$-bounds $(\eta/2, c_R)$, such that
they are $\eta'$-transverse to {\bf 0}. Define
$$ \sigma_{k,t}^j= \sigma_k^j+ t\tau_k^j, ~~ j=1,2.$$
It is easy to check that the sections $\sigma_k^j(t)$ are $\eta/2$-transverse to 
{\bf 0}, therefore its zero sets are isotopic. So we can impose directly that 
$\sigma_k^0\oplus \sigma_k^1$ is transverse to {\bf 0}.
 
Define the function 
\begin{eqnarray*}
F_k: C-Z(\sigma_k^1\oplus \sigma_k^2) & \to & \CP^1 \\
p & \to & \frac{\sigma_k^2(p)}{\sigma_k^1(p)}. 
\end{eqnarray*}
By an analogous argument to that in the precedent lines we assume that $\partial F_k$
is $\eta$-transverse to {\bf 0} away from $Z(\sigma_k^1\oplus \sigma_k^2)$ 
Recall that $F_k^{-1}(0)=Z(\sigma_k^1)$ and $F_k^{-1}(\infty)=
Z(\sigma_k^2)$. Following the proof of Section \ref{proof} 
we obtain a contact fibration $f_k$, for $k$ large enough. But, {\it we do not impose
the local model at a neighborhood of $A$}, therefore we do not perturb any
neighborhood of $A$. We observe that the perturbation at a neighborhood of the set of 
critical points $\Delta$ does not change the counterimages of the points $0$ and $\infty$. 
It follows since there exists a $\rho_0$-neighborhood, being $\rho_0>0$ a uniform constant, 
of these two fibres which
does not intersect the ``bad set'' $\Gamma$, because of the mixed $C^3$-bounds of
the sections $\sigma_k^1$ and $\sigma_k^2$ (which assure contactness at a neighborhood of 
the zero set of the section).
So we perturb the sequence $F_k$ in a $\rho$-neighborhood of $\Delta$, with $0<\rho< \rho_0$.
Therefore, even after performing the needed perturbations, $f_k^{-1}(0)=N_k$ and
$f_k^{-1}(\infty)=N_k'$.

By Proposition \ref{handle} we obtain, following a path between $0$ and $\infty$, that
$f_k^{-1}(0)=N_k$ and $f_k^{-1}(\infty)=N_k'$ are related by a sequence of 
operations of index $n$. This finishes the proof if we admit that the complex structures
coincide for the two sequences of sections.

Suppose now that the sequence of submanifolds $N_k$ is the zero set of a sequence of 
sections $\sigma_k$ with mixed $C^3$-bounds respect to an almost complex structure $J_0$. 
Furthermore $N_k'$ comes from $s_k'$ with mixed $C^3$-bounds respect to $J_1$.
Fix a continuous path $J_t$ joining $J_0$ and $J_1$ in the moduli of compatible
almost complex structures. We use the following
\begin{lemma} \label{small_small}
Let $J_0$ be a compatible almost complex structure in $C$. There exists a uniform 
$\epsilon>0$ satisfying that for any compatible almost 
complex structure $J$ such that $|J-J_0|<\epsilon$, there exist two sequences of sections
$s_k$ and $s_k'$ with mixed $C^3$-bounds, respect to $J_0$ and $J$, which are 
$\eta$-transverse to {\bf 0}, for some $\eta>0$.
Moreover the zero sets of $s_k$ and $s_k'$ are isotopic for $k$ large enough.
\end{lemma}

Find the uniform constants $\epsilon_t=\epsilon$ provided by Lemma \ref{small_small} 
for the almost complex structures $J_t$. This implies, by the continuity of $J_t$, 
that there exists $\epsilon_t'>0$ such that we can find two sequences of sections $s_k^t$ and
$s_k^{t'}$ with mixed $C^3$-bounds respect to $J_t$ and $J_{t+\epsilon_t'}$
whose zero sets are isotopic, for $k$ large enough.

Now we cover the segment $[0,1]$ by open segments $(t, t+\epsilon_t')$. The segment is
compact, so we can find a finite subset $(0, \epsilon_0'), (t_1, t_1+\epsilon_{t_1}'),
\ldots (t_N, 1)$ of the precedent set of segments which covers $[0,1]$.
Obviously, without loss of generality, we can choose $t_j+\epsilon_{t_j}'=t_{j+1}$.

Denote $Z(s)$ the zero set of a given section $s$. Observe that the sets $Z(\sigma_k)$ 
and $Z(s_k^0)$ are related through a cobordism of index $n$, for $k$ large enough. 
But, recall now that $Z(s_k^0)$ and $Z(s_k^{0'})$ are isotopic for $k$ large enough. Again, 
$Z(s_k^{0'})$ and $Z(s_k^{t_1})$ are related through surgeries of index $n$...
Following this argument we find that $Z(\sigma_k^0)$ and $Z(\sigma_k^1)$ are related through
surgeries of index $n$.

Finally recall that if $s_k$ is $\eta$-tranverse to {\bf 0} with mixed $C^r$-bounds
$(c_D, c_R)$ with respect to a complex structure $J$ and to a contact form $\theta$,
then $s_k$ is $\eta/c$-transverse to {\bf 0} with mixed $C^r$-bounds $(cc_D, cc_R)$
with respect to the same contact structure $J$ and to a contact form $\theta'$.
The constant $c>0$ only depends on the contact forms $\theta$ and $\theta'$.
This proves that the result does not depend on the chosen contact form.
\hfill $\Box$

{\bf Proof of Lemma \ref{small_small}:}
Construct a sequence of sections $\sigma_k$ $\eta$-transverse to {\bf 0} using the 
globalization argument of Proposition \ref{global}, for some uniform constant $\eta>0$.
Thus, we obtain
$$ \sigma_k = \sum_{j\in J} w_k \cdot s_{k,x_j}^{ref}, $$
where the points $x_j$ are elements of the set $S$ with the properties described in the 
proof of that Proposition. Recall that the definition of $\sigma_k$ makes sense
if we change the almost complex structure $J_0$ by another one $J$. The obtained
section $\sigma_k'$ will be different because the reference sections $s_{k,x_j}^{ref}$ depend
on the complex structure. However, we claim that $\sigma_k'$ has mixed $C^3$-bounds and
that $\sigma_k-\sigma_k'$ has mixed $C^3$-bounds $(c_u\epsilon, c_R')$,
where $c_u$ and $c_R'$ are uniform constants. Impose
$$ c_u\epsilon < \frac{\eta}{2}. $$
Obviously the constant $\epsilon$ can be chosen uniformly and also $\sigma_k$ and
$\sigma_k'$ have isotopic zero sets. This finishes the proof.

We have only to check that $\sigma_k-\sigma_k'$ have the desired mixed $C^3$-bounds.
We must follow the proof in \cite{IMP99} to check this condition. We only give the
key ideas:
\begin{enumerate}
\item Check that $|s_{k,x, J_0}^{ref}-s_{k,x,J}|_{C^2}=O(\epsilon).$
\item Compute the integer $\CrN$ in the globalization process for the sequence in $J_0$.
Check that this constant is enough also for the globalization process with $J$ if we shrink 
$\epsilon$ uniformly in this globalization process. Why uniformly? Because $\CrN$ is
independent of $k$.
\end{enumerate}
\hfill $\Box$

Proposition \ref{quasi_iso} is weaker than the ``contact Lefschetz hyperplane theorem'' 
proved in \cite{IMP99}, but it is more geometrical and enligthens the behaviour of the
generic sections of the bundles $L^{\otimes k}$. It would be interesting to study the
different connected components of $\CP^1-\Delta$ to control the topology of all the
``approximately holomorphic'' contact zeroes of $L^{\otimes k}$.

\section{The non-exact case.} \label{double}
To conclude the discussion we consider now the non-exact case. The important point is the
following standard result 
\begin{proposition} \label{dou_cov}
Given a non-exact contact manifold $(C,D)$. There exists an exact contact manifold
$(\hat{C},\hat{D})$ which is a non-trivial double covering of $C$. The projection is
a contactomorphism and it can be found a contact form $\hat{\theta}$ in $\hat{C}$ such that
the structure $\Z_2$-action of the covering is a strong anticontactomorphism, i.e.
for $\alpha\in \Z_2$, it is $\alpha_* \hat{\theta}= -\hat{\theta}$.
\end{proposition}
For a simple proof, see \cite{IMP99}. Now, we follow the ideas in \cite{IMP99}. We 
lift $S$ to the double covering and denote it
again by $S$. It is easy to find an almost complex
compatible structure $J$ satisfying that $\alpha_* J=-J$, this implies that $\alpha_*
g_k=g_k$. The application $\alpha$ lifs to a morphism of bundles $\tilde{\alpha}:
L^{\otimes k} \to \bar{L}^{\otimes k}$ (in fact recalling from \cite{IMP99} that $L$
is trivialized by construction, we take the identity in each fiber). This morphism
preserves the connection. So it is easy to check that if $s_k$ is a sequence of
sections of the bundles $S\otimes L^{\otimes k}$ with $C^r$-bounds
$(c_D,c_R)$ with respect to the contact form $\theta$ and the almost complex structure
$J$, then $\alpha_* s_k$ is a sequence of sections of the bundle $S\otimes \bar{L}^
{\otimes k}$ with respect to the contact form $-\theta$ and to the compatible
almost complex structure $-J$ with the same mixed $C^r$-bounds $(c_D,c_R)$.
But, now we identify the bundle $L^{\otimes k}$ with $\bar{L}^{\otimes k}$ using
the anticomplex isomorphism provided by the identity. Then a simple computation
shows that $\alpha_* s_k$ is a sequence os sections of the bundles $S\otimes L^{\otimes k}$
with the same mixed $C^3$-bounds.

Our only task is to assure that all the objects in the construction are $G$-invariant.
The precedent considerations assure this condition in the local constructions.
Now we are going to study the globalization process.
The first point is to achieve the set of points $S$ invariantly. Moreover each $S_j$
has to be $\Z_2$-invariant. This is only true if the action of $\Z_2$ is free, because in any 
other case we would have problems to assure the second property of the set. The perturbation
is performed in a similar way that in the standard case. The key idea is that the 
perturbation term $\tau_{k,x}$
constructed to obtain transversality close to a point $x$ can be transported to the
point $\alpha(x)$ by means of $\tau_{k, \alpha(x)}=\alpha_* \tau_{k,x}$. There is no
interference between $x$ and $\alpha(x)$, because
$$ d_k(x, \alpha(x))=O(k^{1/2}). $$
The term $\tau_{k, \alpha(x)}$ produces the same transversality as $\tau_{k,x}$ because
of the $\Z_2$-invariance of the construction.

The perturbation of the function $F^{s_k}$ can be made invariantly. In the set $A$ it is 
clear. In $\Delta$, the trivialization is invariant, as well as the function $H$ and we 
can choose $w'$ invariant without loss of generality.

So we have constructed a pencil in $\hat{C}$ which is $\Z_2$-invariant. This $\Z_2$-invariance
property allows us to quotient by the group $\Z_2$ obtaining new data
$A=\hat{A}/\Z_2$, $\Delta=\hat{\Delta}/\Z_2$ and $f=f/\Z_2$. It is a trivial exercise
to check that the object so defined is a contact pencil on $C$.
\hfill $\Box$

\end{document}